\documentclass[12pt]{article}

\usepackage[utf8]{inputenc}
\usepackage{amsmath,amsthm,amssymb}

\newcommand{\Homeo}{\mathrm{Homeo\,}}
\newcommand{\uni}{\mathbb{U}}
\newcommand{\E}{\mathrm{E}}
\newcommand{\Aut}{\mathrm{Aut\,}}

\newtheorem{defn}{Definition}[section]
\newtheorem{thm}{Theorem}[section]
\newtheorem{prop}{Proposition}[section]

\begin{document}

\title{Universal spaces for finite group actions on spaces of type $K(\pi,1)$.}
\author{L. Lokutsievskiy}
\maketitle

\begin{abstract}
	In this paper author proposes a construction of a universal space of type $K(\pi,1)$ such that any action (up to homotopy conjugation) of a given finite group $G$ on spaces of the same homotopy type is presented on the constructed space. Moreover, any action of $G$ on any space of type $K(\pi,1)$ is covered by some action of $G$ on this universal space.
\end{abstract}

\section{Introduction}

First of all let us consider some key notions.

\begin{defn}
	Let $p_1$ and $p_2$ be actions\footnote{The term ``group action'' is used in standard sense. It is a homomorphism $p:G\to\Homeo(T)$ where $G$ is a group and $T$ is a topological space.} of an abstract group $G$ on topological spaces $T_1$ and $T_2$ respectively. We say that $p_1$ covers $p_2$ if there exist homotopy equivalence $\varphi:T_1\to T_2$ and an automorphism $\theta\in\Aut G$ such that
	
	$$
		\forall g\in G\qquad \varphi\circ p_1(g) = p_2\bigl(\theta(g)\bigr)\circ \varphi.
	$$
\end{defn}

Obviously, the binary relation of covering is reflexive and transitive. But in general it is not symmetric and it is easy to construct counter-examples even if $T_1$ and $T_2$ have CW-complex structure.

\begin{defn}
	We say that the actions $p_1$ and $p_2$ are homotopy conjugate if $p_1$ and $p_2$ cover each other.
\end{defn}

The following main theorem is proved in this article:

\begin{thm}
\label{thm_main}
	Let $\pi$ be an abstract group and $G$ be a finite group. Then there exists a CW-complex $\uni(G,\pi)$ of type $K(\pi,1)$ such that for any CW-complex $T$ of the same homotopy type and any regular\footnote{The theorem is proved with minimal natural restrictions like CW-complex structure. See definitions \ref{defn_free_regular} and \ref{defn_nonfree_regular} for more details} action $p$ of $G$ on $T$ there exists a free action $p_\uni$ of $G$ on $\uni(G,\pi)$ such that $p$ is covered by $p_\uni$. Moreover, if $p$ is free then $p_\uni$ and $p$ are homotopy conjugate.
\end{thm}

Let us also remark that there exists a theorem of full classification of regular free action of $G$ on spaces of type $K(\pi,1)$ (see \cite{Lokut_1}). This classification is closely connected to extensions of $G$ and its cohomology (see \cite{Brown}, IV.6)

\section{Classification of free actions on spaces of type $K(\pi,1)$}

The proof of theorem \ref{thm_main} is based on the classification of regular free actions (for the proofs of the results described in this paragraph, see \cite{Lokut_1}). 

\begin{defn}
\label{defn_free_regular}
	Let $p$ be a free action of a finite group of $G$ on a CW-complex $T$. Then $p$ is called regular if the quotient space $T/p(G)$ is also a CW-complex\footnote{In fact, since the projection $T\to T/p(G)$ is a regular covering, we see that any CW-complex structure on $T/p(G)$ induces a CW-complex structure on the space $T$.}.
\end{defn}

\begin{prop}
	The binary relation of covering for regular free actions on spaces of type $K(\pi,1)$ is symmetric. Consequently in this case it is an equivalence relation and it coincides with homotopy conjugation.
\end{prop}

\begin{thm}
\label{thm_classification}
	The set of all regular free actions (up to homotopy conjugation) of $G$ on spaces of type $K(\pi,1)$ is in one-to-one correspondence with the set of all classes of equivalent extensions of $G$ by $\pi$. The following extension (the final part of the long exact sequence of the covering $T\to T/p(G)$) corresponds to the action $p$:

	$$
		1\to \pi \to \pi_1\bigl(T/p(G)\bigr) \to G \to 1.
	$$
\end{thm}

It  is important to notice that two extensions of $G$ by $\pi$ are equivalent if they are isomorphic in the category of short exact sequences.

Full classification of group extensions was constructed by Eilenberg and MacLane in 1947 (see \cite{EilenbergMacLane}). It is closely connected to cohomology of $G$ and can be obviously used together with theorem \ref{thm_classification}\footnote{Eilenberg and MacLane worked with classes of congruent extensions. Classes of equivalent extensions is used in theorem \ref{thm_classification}. However if two extensions are congruent then they are isomorphic. So the results of Eilenberg and MacLane can be used with theorem \ref{thm_classification}.}.

Before we start the proof of the main theorem \ref{thm_main} let us make a definition of regular non-free action.

\begin{defn}
\label{defn_nonfree_regular}
	A non-free action of $G$ on $T$ is called regular if there exists a contractible space $\E G$ with a free action of $G$ on it such that the diagonal (free) action of $G$ on $T\times\E G$ is regular.
\end{defn}

\section{The proof of theorem \ref{thm_main}}

The proof of theorem \ref{thm_main} is not complicated and its main idea is the following: we will construct a space $\uni(G,\pi)$ with all actions of $G$ on it (using the classification theorem \ref{thm_classification}) instead of searching the necessary covering action $p_\uni$ on $\uni(G,\pi)$ for each CW-complex of type $K(\pi,1)$ and any action $p$ on it.

The proof is in 3 steps:

\begin{enumerate}
	\item For any extension $e$ of $G$ by $\pi$ we will construct a CW-complex $\uni(G,\pi,e)$ of type $K(\pi,1)$ and a free action $p_e$ of $G$ on it such that $e$ corresponds to $p_e$ by theorem \ref{thm_classification}.
	\item We will show that for any two extensions $e_1$ and $e_2$ the spaces $\uni(G,\pi,e_1)$ and $\uni(G,\pi,e_2)$ are homeomorphic.
	\item For any space $T$ and action $p$ on it we will construct the required covering action $p_\uni$.
\end{enumerate}

\begin{proof}

{\bf STEP 1.} Consider an extension $e:1\to\pi\to S\to G\to 1$. Let us use the Milnor construction (for example see \cite{Husemoller}, IV.11): put $\E S=S\ast S\ldots$ ($\E S$ is a join of a countable number of $S$ where $S$ is considered as a topological space with discrete topology). Then $\E S$ is a contractible space and $S$ acts freely and continuously on it. Trivially, $\pi\subseteq S$ and the projection $\E S\to \E S/\pi$ is a regular covering. By definition, put $\uni(G,\pi,e)=\E S/\pi$. Obviously, the space $\uni(G,\pi,e)\simeq K(\pi,1)$ is a CW-complex\footnote{Moreover, $\uni(G,\pi,e)$ is a simplicial complex. See \cite{Bredon}, III.1 for more details.}. Since $S$ acts freely on $\E S$, we see that the induced action $p_e$ of $G$ on $\uni(G,\pi,e)$ is free and continuous\footnote{The action $p_e$ is continuous since the projection $\E S \to \uni(G,\pi,e)$ is an open map.}. Moreover, $\uni(G,\pi,e)/G=\E S/S$ is a CW-complex. Hence $p_e$ is regular.

Since  $G$ is finite, we see that the projection $\uni(G,\pi,e)\to \E S/S$ is a regular covering. Let us write out the long exact sequence of this covering:

$$
	\ldots\to 1\to \pi\to S\to G\to 1.
$$

So $e$ corresponds to the action $p_e$ on $\uni(G,\pi,e)$ by theorem \ref{thm_classification}. Q.E.D.

{\bf STEP 2.} Let us consider two extensions $e_i:1\to\pi\to S_i\to G\to 1$, $i=1,2$ and show that the spaces $\uni(G,\pi,e_1)$ and $\uni(G,\pi,e_2)$ are homeomorphic.

The group $\pi$ acts freely on $S_1$ and $S_2$ by left shifts. Moreover, the sets $S_i/\pi$ of orbits and the set $G$ are of equal cardinality. Hence there exists a set-theoretic bijection $f:S_1\to S_2$ and we can assume that $f$ respects the actions of $\pi$ on $S_1$ and $S_2$. The bijection $f$ induces a homeomorphism $\widetilde f:\E S_1\to \E S_2$ and $\widetilde f$ respects the actions of $\pi$ on these spaces. Consequently the map $\widetilde f$ provides a homeomorphism $\uni(G,\pi,e_1)\to\uni(G,\pi,e_2)$. Q.E.D.

By definition, put $\uni(G,\pi)=\uni(G,\pi,e)$ where $e$ is an arbitrary extension of $G$ by $\pi$ (for example, $e$ is the extension of direct product).

So we proved that for any extension $e$ there exists an action $p_e$ of $G$ on $\uni(G,\pi)$ such that $e$ corresponds to $p_e$ by theorem \ref{thm_classification}

{\bf STEP 3.} Consider now a regular non-free action $p$ of $G$ on a CW-complex $T$ of type $K(\pi,1)$. Let us use the Borel construction: there exists an extension $e$ of $G$ by $\pi$ corresponding to the free diagonal action $p'$ of $G$ on $T\times\E G$.

It follows from theorem \ref{thm_classification} that $p'$ and $p_e$ are homotopy conjugate. Hence there exists a homotopy equivalence $\psi:\uni(G,\pi)\to T\times\E G$ and $\psi$ respects the actions of $G$ up to an automorphism $\theta\in\Aut G$. It is easy to see that $p_e$ covers $p$ and the composition of $\psi$ and the projection $T\times\E G\to T$ is the necessary covering homotopy equivalence.

If the action $p$ is free then the construction of $p'$ is not needed since $p$ and $p_e$ are already homotopy conjugate.

\end{proof}

I was asked about existence of universal spaces for finite group actions during my talk on the Postnikov research seminar "Algebraic topology and its applications" directed by V.M. Buhshtaber. This article is an answer to this question and I am very grateful to participants of the seminar for this question.


\begin{thebibliography}{99}

\bibitem{Brown} {\it K.S. Brown}, ``Cohomolgy of groups'', Springer-Vergal, New York, Heidelberg, Berlin, 1982.

\bibitem{Bredon} {\it G.E. Bredon}, ``Introduction to compact transformation group'', Academic press, New York -- London, 1972

\bibitem{Husemoller} {\it D. Huemoller}, ``Fibre Bundle'', McGraw-Hill Book Company, New York etc., 1966

\bibitem{Lokut_1} {\it Lev Lokutsievskiy}, ``Homotopy classification of finite group actions on aspherical spaces'', 2010, arXiv:1008.2980 (will be printed in ``Journal of Mathematical Science'').

\bibitem{EilenbergMacLane} {\it S. Eilenberg, S. MacLane }, ``Cohomology theory in abstract groups. II. Group extensions with a non-abelian kernel'', Ann. Math., 1947, (2) 48, p. 199-236.

\end{thebibliography}
\end{document}